\documentclass[12pt]{amsart}
\usepackage{geometry}                
\usepackage{moreverb,url}
\usepackage{amssymb,amsmath,amsfonts}
\usepackage{bm}
\usepackage[font=scriptsize, labelfont=bf]{caption}
\usepackage{float}
\usepackage{afterpage}
\usepackage{subcaption}
\usepackage[squaren]{SIunits}
\usepackage{verbatim}
\usepackage[square,numbers]{natbib}
\usepackage{stmaryrd}
\usepackage{fancyhdr}
\usepackage{syntax,etoolbox}
\usepackage{graphicx}
\usepackage{hyperref}
\hypersetup{colorlinks=true,allcolors=blue}
\usepackage[square,numbers]{natbib}
\usepackage{lineno}

\begin{document}

	\title[Virus Deformation]{Modelling Virus Contact Mechanics under Atomic Force Imaging Conditions}

	\author{Paolo Piersanti}
	\address{Department of Mathematics and Institute for Scientific Computing and Applied Mathematics, Indiana University Bloomington, 729 East Third Street, Bloomington, Indiana, USA}
	\email{ppiersan@iu.edu}

	\author{Kristen White}
	\address{Department of Chemistry, Indiana University Bloomington, 800 E Kirkwood Ave, Bloomington, Indiana, USA}
	\email{kw98@iu.edu}
	
	\author{Bogdan Dragnea}
	\address{Department of Chemistry, Indiana University Bloomington, 800 E Kirkwood Ave, Bloomington, Indiana 47405, USA}
	\email{dragnea@iu.edu}
		
	\author{Roger Temam}
	\address{Department of Mathematics and Institute for Scientific Computing and Applied Mathematics, Indiana University Bloomington, 729 East Third Street, Bloomington, Indiana, USA}
	\email[Corresponding author]{temam@indiana.edu}

\begin{abstract}
In this paper we present a discrete model governing the deformation of a convex regular polygon subjected not to cross a given flat rigid surface, on which it initially lies in correspondence of one point only. First, we set up the model in the form of a set of variational inequalities posed over a non-empty, closed and convex subset of a suitable Euclidean space. Secondly, we show the existence and uniqueness of the solution. The model provides a simplified illustration of processes involved in virus imaging by atomic force microscopy:  adhesion to a surface, distributed strain, relaxation to a shape that balances adhesion and elastic forces. The analysis of numerical simulations results based on this model opens a new way of estimating the contact area and elastic parameters in virus contact mechanics studies.
\end{abstract}

\today

\maketitle

\section{Introduction}

Due to its relevance for a variety of natural phenomena, contact mechanics of nanoscopic, isometric polyhedral shells is a topic of great interest. For instance, formation of polyhedral cages of water molecules encapsulating other chemical species at low temperatures and high pressures govern the sequestration of hydrocarbon gas molecules on arctic sea floors, and in permafrost \cite{Sloan2003}. Bulk polycrystalline samples of such gas hydrates (clathrates) can be 20-40 times stronger than ice under uniaxial compression. At larger spatial scales, but still nanoscopic, highly symmetric virus protein cages encapsulate nucleic acid cargo, exhibiting a hardness that is responsive to the chemical environment. For such examples the discrete nature of the subunits, and the size-dependence of their properties cannot be ignored when attempting to understand cage mechanics in response to chemical or physical forces. 

During the virus life cycle, there are multiple instances of interactions with interfaces. The ensuing mechanical stresses elicit a response from either the host or the virus \cite{Greber2014}. Since the nature and magnitude of resultant deformations are key to understanding the response, virus mechanics has been the topic of both experimental and computational studies \cite{Guerra2017, Zeng2017a, Zeng2021}. In many viruses, the attractive interactions between the protein building blocks of a virus shell (also known as the capsid) are necessarily weak, at least at the assembly stage \cite{Prevelige:1993db, Zlotnick:1999ww}. This allows for error correction, during the growth phase. By comparison, the oligomeric building blocks themselves tend to be much more stable against dissociation \cite{Dragnea2019,Zlotnick:2007ji}. As a consequence of the difference in intra- and inter-subunit interaction potential energies, and of the geometric frustration, local strain is expected to vary across the shell structure \cite{Aggarwal2012, Zandi2005}. The inhomogeneous elastic stress distribution  will influence how capsids deform under mechanical stress.

Early theories of mechanical deformation in viruses sought to explain observations from nanoindentation experiments \cite{Michel:2006jt, Ivanovska2004}. Nanoindentation experiments involve the measurement of axial compression as a function of the force exerted by an atomic force microscope (AFM) probe in a direction normal to the substrate on which the virus sits \cite{Roos2011, Thompson2020}. Early models were based on the continuum theory of elasticity \cite{Gibbons2007, Bruinsma2015}. Later on, the discrete nature of the building blocks was taken into account in coarse-grained and all-atom molecular dynamics models of nanoindentation \cite{Vliegenthart2006, Arkhipov:2009kt, Mannige2009, Krishnamani2016}. These models provided better qualitative predictions for stiffness distribution across the capsid than the continuum theory, at least for small viruses. They are also inherently capable of including fluctuations, which can influence the path through the energy landscape as the virus is deformed. However, the great complexity of some of the numerical models make it difficult to evidence the layers that would lead to a conceptual understanding of the process. The challenges are even greater if one considers that, in nanoindentation experiments, the read-out of a three-dimensional conformational change is one-dimensional (axial distance). Solving the inverse problem of finding material mechanical parameters from the axial distance dependence depends on the model. Moreover, since the exact AFM probe-virus contact geometry changes from particle to particle and probe to probe, force-deformation curves usually show significant experimental spread.

In this paper we propose a new way to perform a three-dimensional measurement of conformational changes of a virus under directional pressure, and extract elastic parameters via an elastic model based on a system of coupled discrete subunits with localized interactions. 
The model we will be considering is discrete, and is governed by a set of variational inequalities posed over a non-empty, closed, and convex subset of the Euclidean space.
Similar problems have been considered in the continuum framework in the static case~\cite{CiaPie2018b,CiaPie2018bCR,CiaMarPie2018b,CiaMarPie2018,Pie2020-1,PS} and in the time-dependent case~\cite{Pie2020}.
For clarity and simplicity, we consider the geometry two-dimensional, and the interactions linear. Extensions to three dimensions and non-linear interactions will be considered in the forthcoming paper~\cite{PWDT3D}.

\emph{This article is dedicated to Robert P. Gilbert on the occasion of his 90th birthday, with friendship and much appreciation for his scientific contributions and for his services to the community.}

\section{Basic principles of atomic force microscopy}
 
\begin{figure}
\includegraphics[width = \textwidth]{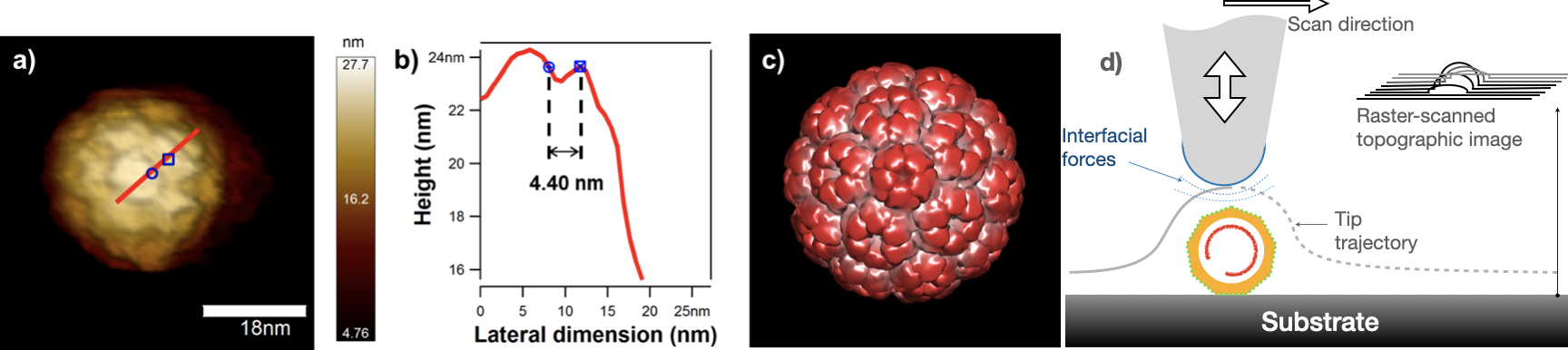}
\caption{a) Example of an AFM image of a BMV particle in aqueous medium. b) Cross-section showing the shape profile and typical attainable resolution. c) Molecular model (reprinted with permission from ref.~\cite{Zeng2017a}. d) Schematic of the atomic force probe-sample geometry. The probe is in the shape of a rounded, inverted cone (gray). The sample is formed of a virus (orange) adsorbed on the surface of a flat, solid substrate. As the probe tip is laterally scanned over the substrate, the probe-sample distance is kept constant, and the tip trajectory relative to a reference position is recorded.}
\label{afm}
\end{figure}

Atomic force microscopy (AFM) is currently the only method capable of imaging small virus morphology in real-time under environmental conditions, with a spatial resolution routinely reaching 3-5 nm \cite{Kuznetsov2011, Baclayon2010, Cartagena2013, Zeng2017a, Calo2021}. The study of single virus particles by atomic force microscopy and force spectroscopy has led to insights in virus mechanics and its relationship with the virus structure and the chemical environment \cite{Mateu2012, Castellanos2012a, Pang2013a,Snijder2013a,Kononova2013,Guerra2017,Zeng2021}. A schematic of the atomic force imaging principle is presented in Fig.~\ref{afm}.

In AFM, a sharp (several nm) probe is brought in the proximity of the sample surface until the probe-sample distance is so small that interfacial forces between the two become measurable. By raster scanning the probe over the sample while adjusting the probe-sample distance according the sample topography to maintain a constant interfacial force, one generates a nanoscopic map of the sample topography. Scan areas may span from several $\mathrm{nm^2}$ to several hundreds of $\mathrm{\mu m^2}$. 

In a direction normal to the substrate, spatial resolution is typically better than a nm. However, lateral resolution is limited by the sample topographic range and the probe radius, Fig~\ref{afm}.

Since the normal pressure can be adjusted at will during imaging, one can create 3-dimensional topographic maps of the upper shell under several constant forces. Our hypothesis is that, from these topographic maps one should be able to determine the mechanical parameters that describe bending and stretching of the shell under deformation.

\section{Geometrical preliminaries}
\label{Sec:1}

\begin{figure}[ht]
	\includegraphics[width=\linewidth]{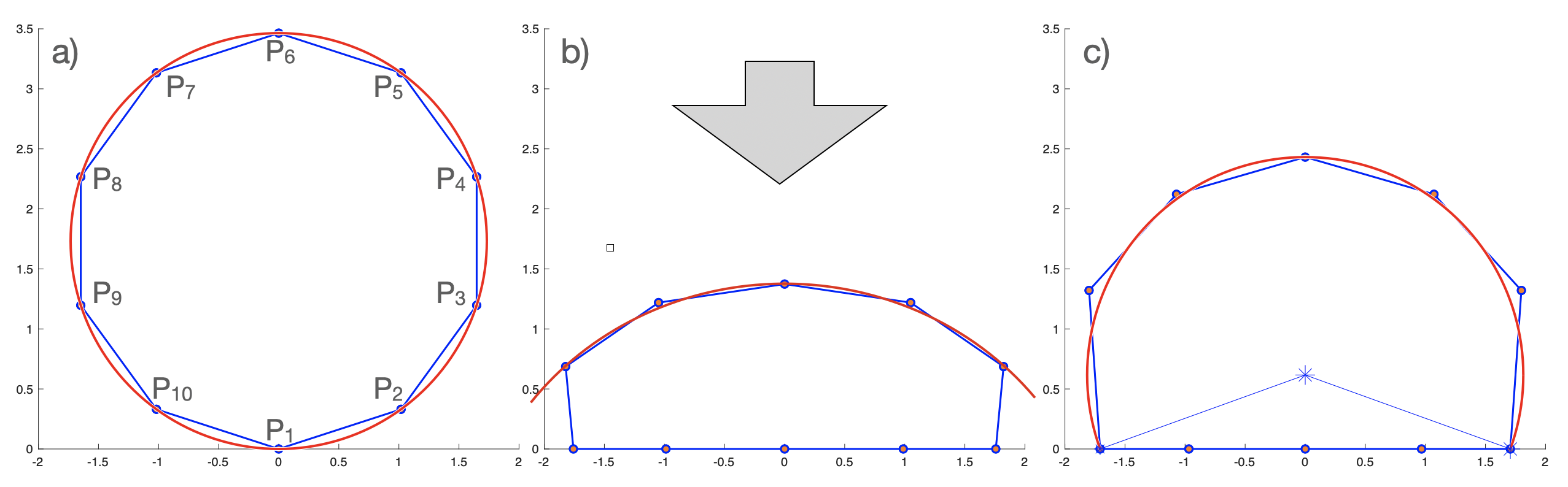}
	\caption{a) A convex regular polygon with $n=10$ edges initially in contact with the surface at one fixed vertex. b) The compressed shape under a constant, uniaxial (normal) force. Red: the circular sector that best fits the free perimeter. c) The final shape after removal of the normal force under the constraint of irreversible adhesion. }
	\label{geometry}
\end{figure}

In this paper we consider the deformation of a convex regular polygon with $n$ edges whose vertices (and so the edges) are subjected not to cross an \emph{undeformable} flat surface. We assume that one and only one vertex of the undeformed reference configuration is initially in contact with the flat surface, and we denote this vertex by $P_1$. The vertices of the polygon are labelled counter-clockwise from $P_1$ to $P_n$, Fig.~\ref{geometry}a.
We further \emph{assume} that the vertex $P_1$ undergoes no displacement; this assumption is critical to establish the existence and uniqueness of the solution of the governing equations. We minimize the energy associated with the deformation of the convex regular polygon with $n$ edges under two types of constraints: i) under the action of an applied body force, Fig.~\ref{geometry}b, and ii) after the force is removed, under the assumption of irreversible adhesion, i.e. vertexes brought on the rigid surface by the pressure in i) will remain on the surface at force removal, Fig.~\ref{geometry}c. We consider the elastic energy to have a quadratic dependence on the elongation of the edges as well as the variation of the angle between any pairs of consecutive edges~\cite{Zeng2017a}. Furthermore, we assume that the edges are \emph{massless} and can \emph{only stretch}, in the sense that \emph{there is no torsion acting on them}.

A point $A$ in the plane corresponds to a vector in $\mathbb{R}^2$ of the form $A=\begin{pmatrix}x_A\\y_A\end{pmatrix}$, where $x_A \in \mathbb{R}$ denotes the abscissa of the point and $y_A \in \mathbb{R}$ denotes its ordinate.
Let us consider a Cartesian frame for the two-dimensional plane with origin $O=\begin{pmatrix}0\\0\end{pmatrix}$ and with canonical directions $\vec{e}_1=\begin{pmatrix}1\\0\end{pmatrix}$ and $\vec{e}_2=\begin{pmatrix}0\\1\end{pmatrix}$.

The position vector associated with the point $A$ is denoted by $\overrightarrow{OA}$; the transpose of a vector $\overrightarrow{OA}$ is denoted by $\overrightarrow{OA}^T=(x_A,y_A)$; the angle between three points $A$, $B$ and $C$ with vertex at $B$ is either denoted by $\widehat{ABC}$ or by a Greek letter.
The Euclidean inner product and the vector product between two vectors $\overrightarrow{OA}$ and $\overrightarrow{OB}$ are respectively denoted by $\overrightarrow{OA} \cdot \overrightarrow{OB}=\overrightarrow{OA}^T \overrightarrow{OB}$ and $\overrightarrow{OA} \times \overrightarrow{OB}$. The Euclidean norm of $\overrightarrow{OA}$ is denoted $\left|\overrightarrow{OA}\right|$. Matrices, apart from the identity matrix $\mathcal{I}$ and the square null matrix $\mathcal{O}$, are denoted by capital Greek letters.
Tensors are denoted by boldface capital Latin letters.

Let $n \ge 3$ be an integer number. A regular polygon with $n$ edges is the portion of plane within a non-self intersecting closed broken line, whose $n$ edges have all the same length (cf., e.g., \cite{Kiselev2006}).

\section{Formulation and well-posedness of the two-dimensional discrete model}
\label{Sec:2}

When a vertex $P_i$ of the convex regular polygon with $n$ edges under consideration undergoes the action of an applied body force, we denote by $P_i'$ the new coordinates of the vertex in the plane. 

Let $\bm{F}=(\vec{f}_i)_{i=2}^n \in \mathbb{R}^{2n-2}$ denote the array of applied body forces acting on the polygon, where the vector $\vec{f}_i$ denotes the applied body force acting on the vertex $P_i$. The application of the force vector $\vec{f}_i$ on the point $P_i$ displaces the position vector $\overrightarrow{OP_i}$ by a vector $\vec{u}_i$, and transforms the  vector $\overrightarrow{OP_i}$ into the  vector $\overrightarrow{OP_i'}$ via the following relation:
\begin{equation*}
\overrightarrow{OP_i'} = \overrightarrow{OP_i} + \vec{u}_i,\quad\textup{ for each } 2 \le i \le n.
\end{equation*}

We denote by $\ell$ the length of any edge of the undeformed reference configuration of the convex regular polygon with $n$ edges under consideration, namely,
\begin{equation*}
\left|\overrightarrow{P_i P_{i+1}}\right| =\ell,
\end{equation*}
where the indices are meant from now on \emph{modulo} $n$.

Since the point $P_1$ undergoes, by assumption, no deformation we let $\vec{u}_1=(0,0)^T$.
The stretching energy associated with the displacement 
$$
\bm{U}=\begin{pmatrix}
\vec{u}_2\\
\vdots\\
\vec{u}_n
\end{pmatrix} \in \mathbb{R}^{2n-2},
$$ 
is computed via Hooke's law (cf., e.g., \cite{TemamMiranville2005}), i.e.,
\begin{equation*}
J_s(\bm{U}):=\dfrac{k}{4} \sum_{i=1}^{n} \sum_{j \in \{i-1,i+1\}} |\vec{u}_i - \vec{u}_j|^2,
\end{equation*}
where the elastic constant $k>0$ is associated with the elongation properties of the constitutive material, and the nature of the energy is aptly recalled by the subscript ``$s$''.

The stretching energy $J_s(\bm{U})$ can equivalently be expressed in matrix form. The matrix associated with the stretching elastic force is a $2n \times (2n-2)$ matrix
\begin{equation*}
\Sigma=
\begin{pmatrix}
-\mathcal{I} & \mathcal{O} & \mathcal{O} & \dots & \mathcal{O} \\
\mathcal{I} & -\mathcal{I} & \mathcal{O} & \dots & \mathcal{O} \\
\vdots &  &  & \ddots & \\
\mathcal{O}& \dots&\mathcal{O}&\mathcal{I}&-\mathcal{I}\\
\mathcal{O}& \dots&\mathcal{O}&\mathcal{O}&\mathcal{I}
\end{pmatrix}
.
\end{equation*}

Therefore, the stretching energy $J_s(\bm{U})$ in matrix form is
\begin{equation*}
J_s(\bm{U})=\dfrac{k}{2} \bm{U}^T \Sigma^T \Sigma \bm{U}.
\end{equation*}

The matrix $\Sigma^T \Sigma$ is positive-definite, in the sense that the smallest eigenvalue is greater than zero. To see this, it suffices to observe that the matrix $\Sigma$ has the same structure as the discrete gradient matrix (cf., e.g., \cite{Quarteroni2010}).

The variation of the angle between two consecutive edges is also associated with a change in the energy. 
From now one, we will refer to this kind of energy as \emph{bending energy}.
Denote by $\alpha_i=\widehat{P_{i-1} P_i P_{i+1}}$ the angle between the points $P_{i-1}$, $P_i$ and $P_{i+1}$ that intersects the interior of the polygon under consideration.

If the action of an applied body forces changes the angle $\alpha_i$ into the angle $\alpha_i'=\widehat{P_{i-1}' P_i' P_{i+1}'}$, the corresponding bending energy is given by
\begin{equation*}
\dfrac{\kappa}{2} |\alpha_i-\alpha_i'|^2,
\end{equation*}
where the elastic constant $\kappa>0$ is associated with the bending properties of the constitutive material.


If the difference between $\alpha_i$ and $\alpha_i'$ is small, we can approximate $\alpha_i-\alpha_i' \approx \sin(\alpha_i - \alpha_i')$. The latter has the advantage that it can be expressed in terms of a vector product. More specifically, we have:
\begin{equation*}
\begin{aligned}
\alpha_i-\alpha_i' & \approx \sin(\alpha_i -\alpha_i')=\sin \alpha_i \cos \alpha_i'-\cos\alpha_i \sin\alpha_i'\\
&=\dfrac{|\overrightarrow{P_iP_{i+1}} \times \overrightarrow{P_{i-1}P_i}|}{|\overrightarrow{P_iP_{i+1}}| |\overrightarrow{P_{i-1}P_i}|}
\dfrac{\overrightarrow{P_{i-1}'P_i'} \cdot \overrightarrow{P_i'P_{i+1}'}}{|\overrightarrow{P_i' P_{i+1}'}| |\overrightarrow{P_{i-1}'P_i'}|}
-\dfrac{\overrightarrow{P_{i-1}P_i} \cdot \overrightarrow{P_iP_{i+1}}}{|\overrightarrow{P_i P_{i+1}}| |\overrightarrow{P_{i-1}P_i}|}
\dfrac{|\overrightarrow{P_i'P_{i+1}'} \times \overrightarrow{P_{i-1}'P_i'}|}{|\overrightarrow{P_i'P_{i+1}'}| |\overrightarrow{P_{i-1}'P_i'}|}.
\end{aligned}
\end{equation*}

Let us observe that, if the displacement $\bm{U}=(\vec{u}_i)_{i=2}^n$ is infinitesimal, we have
\begin{equation*}
|\overrightarrow{P_{i-1}'P_i'}| \approx |\overrightarrow{P_{i-1}P_i}|.
\end{equation*}

Since
$$
\overrightarrow{P_{i-1}'P_i'} = \overrightarrow{OP_i'} - \overrightarrow{OP_{i-1}'}=\overrightarrow{OP_i}+\vec{u}_i-\overrightarrow{OP_{i-1}}-\vec{u}_{i-1}
=\overrightarrow{P_{i-1}P_i}+ \vec{u}_i -\vec{u}_{i-1},
$$
then the properties of the vector product (cf., e.g., \cite{Strang1980}) in turn imply:
\begin{equation*}
|\overrightarrow{P_i'P_{i+1}'} \times \overrightarrow{P_{i-1}'P_i'}|
=|(\overrightarrow{P_iP_{i+1}}+ \vec{u}_{i+1} -\vec{u}_i) \times (\overrightarrow{P_{i-1}P_i}+ \vec{u}_i -\vec{u}_{i-1})|\approx
|\overrightarrow{P_iP_{i+1}} \times \overrightarrow{P_{i-1}P_i}|.
\end{equation*}

Therefore, the angle variation reads
\begin{equation*}
\begin{aligned}
\alpha_i-\alpha_i' &\approx \dfrac{|\overrightarrow{P_iP_{i+1}} \times \overrightarrow{P_{i-1}P_i}|}{|\overrightarrow{P_iP_{i+1}}|^2 |\overrightarrow{P_{i-1}P_i}|^2}
\left(\overrightarrow{P_{i-1}'P_i'} \cdot \overrightarrow{P_i'P_{i+1}'}-\overrightarrow{P_{i-1}P_i} \cdot \overrightarrow{P_iP_{i+1}}\right)\\
&=\dfrac{|\overrightarrow{P_iP_{i+1}} \times \overrightarrow{P_{i-1}P_i}|}{|\overrightarrow{P_iP_{i+1}}|^2 |\overrightarrow{P_{i-1}P_i}|^2}
\left((\overrightarrow{P_{i-1}P_i}+ \vec{u}_i -\vec{u}_{i-1}) \cdot (\overrightarrow{P_iP_{i+1}}+ \vec{u}_{i+1} -\vec{u}_i)-\overrightarrow{P_{i-1}P_i} \cdot \overrightarrow{P_iP_{i+1}}\right)\\
&\approx\dfrac{|\overrightarrow{P_iP_{i+1}} \times \overrightarrow{P_{i-1}P_i}|}{|\overrightarrow{P_iP_{i+1}}|^2 |\overrightarrow{P_{i-1}P_i}|^2}
\left(\overrightarrow{P_{i-1}P_i} \cdot (\vec{u}_{i+1} -\vec{u}_i)+ \overrightarrow{P_iP_{i+1}} \cdot(\vec{u}_i -\vec{u}_{i-1})\right)\\
&=\dfrac{|\overrightarrow{P_iP_{i+1}} \times \overrightarrow{P_{i-1}P_i}|}{|\overrightarrow{P_iP_{i+1}}|^2 |\overrightarrow{P_{i-1}P_i}|^2}
\left(\overrightarrow{P_{i-1}P_i} \cdot \vec{u}_{i+1}+(\overrightarrow{P_iP_{i+1}}-\overrightarrow{P_{i-1}P_i})\cdot \vec{u}_i-\overrightarrow{P_iP_{i+1}} \cdot\vec{u}_{i-1}\right).
\end{aligned}
\end{equation*}

Letting $C:=\frac{|\overrightarrow{P_iP_{i+1}} \times \overrightarrow{P_{i-1}P_i}|}{|\overrightarrow{P_iP_{i+1}}|^2 |\overrightarrow{P_{i-1}P_i}|^2}$ and observing that this quantity is independent of the index $i$, we can express the total bending energy in terms of the vertices displacements:
\begin{equation*}
J_b(\bm{U})=\dfrac{\kappa C^2}{2} \sum_{i=2}^n \left|\overrightarrow{P_{i-1}P_i} \cdot \vec{u}_{i+1}+(\overrightarrow{P_iP_{i+1}}-\overrightarrow{P_{i-1}P_i})\cdot \vec{u}_i-\overrightarrow{P_iP_{i+1}} \cdot\vec{u}_{i-1}\right|^2.
\end{equation*}

The bending energy $J_b(\bm{U})$ can equivalently be expressed in matrix form. The matrix associated with the bending elastic force is a $n \times (2n-2)$ matrix:
\begin{equation*}
	\hspace{-1.5cm}\Theta=
	\begin{pmatrix}
		(\overrightarrow{P_2P_3}-\overrightarrow{P_1P_2})^T & \overrightarrow{P_1P_2}^T & 0 & \dots & 0 \\
		\overrightarrow{P_3P_4}^T & (\overrightarrow{P_3P_4}-\overrightarrow{P_2P_3})^T & \overrightarrow{P_2P_3}^T & \dots & 0 \\
		\vdots &  &  & \ddots & \\
		0& \dots&0&(\overrightarrow{P_{n-1}P_n}-\overrightarrow{P_{n-2}P_{n-1}})^T&\overrightarrow{P_{n-2}P_{n-1}}^T\\
		0& \dots&0&\overrightarrow{P_nP_1}^T&(\overrightarrow{P_nP_1}-\overrightarrow{P_{n-1}P_n})^T
	\end{pmatrix}
	.
\end{equation*}

Therefore, the bending energy $J_b(\bm{U})$ in matrix form is
\begin{equation*}
	J_b(\bm{U})=\dfrac{\kappa C^2}{2} \bm{U}^T \Theta^T \Theta \bm{U},
\end{equation*}
where the nature of the energy is aptly recalled by the subscript ``$b$''.

The matrix $\Theta^T \Theta$ is nonnegative-definite, in the sense that the smallest eigenvalue is greater or equal than zero.
The total elastic energy associated with the displacement tensor $\bm{U}$ thus takes the following form
\begin{equation*}
J(\bm{U})=J_s(\bm{U}) +J_b(\bm{U})= \bm{U}^T \left(\dfrac{k}{2} \Sigma^T \Sigma+\dfrac{\kappa C^2}{2} \Theta^T \Theta\right) \bm{U},
\end{equation*}
so that the total energy matrix is positive-definite, being $\Sigma^T \Sigma$ positive-definite. This makes the total elastic energy a strictly convex quadratic functional.

The search for an equilibrium position for the deformed polygon amounts to minimizing the corresponding total elastic energy.
In view of the geometrical constraint according to which the vertices do not have to cross the given flat surface, the admissible displacement fields are to be sought in the following set
\begin{equation*}
\mathcal{U}:=\left\{\bm{V}=(\vec{v}_i)_{i=2}^n \in \mathbb{R}^{2n-2}; \vec{v}_i=\begin{pmatrix}v_{i,1}\\v_{i,2}\end{pmatrix} \textup{ and }(\overrightarrow{OP_i}+\vec{v}_i)\cdot \vec{e}_2 \ge 0 \textup{ for all }2 \le i \le n\right\}.
\end{equation*}

It is straightforward to observe that the set $\mathcal{U}$ is non-empty (as $\bm{V}=\bm{0}\in \mathbb{R}^{2n-2}\in \mathcal{U}$), closed and convex.

Therefore, the latter together with the fact that the total elastic energy functional is strictly convex, imply that the quadratic minimization problem 
\begin{equation*}
\inf_{\bm{V} \in \mathcal{U}} (J(\bm{V}) -\bm{F}^T \bm{V})
\end{equation*}
admits a unique minimizer (cf., e.g., \cite{EkelandTemam1999}). Finding the solution for this minimization problem is equivalent to finding a tensor $\bm{U}$ that solves the following variational inequalities~\cite{DuvLions}:
\begin{equation*}
\left(\left(k \Sigma^T \Sigma+\kappa C^2 \Theta^T \Theta\right) \bm{U}\right) \cdot (\bm{V}-\bm{U}) \ge \bm{F} \cdot (\bm{V}-\bm{U}),
\quad\textup{ for all } \bm{V} \in \mathcal{U}.
\end{equation*}

\section{Numerical experiments. Part~{I}: Analysis of the deformation under the action of an applied body force}
\label{Sec:3}

The first numerical experiment we conduct on the proposed model is classical, and amounts to finding the position of the deformed reference configuration of the convex regular polygon undergoing the action of an applied body force which acts on each vertex with the same magnitude. This distributed applied force models the downward pressure that an AFM probe, which is comparable in size with the virus will exert.

The numerical implementation for the problem under consideration is carried out by resorting to the primal-dual active set method (cf., e.g., \cite{SunYuan2006}). We consider different instances of the array $\bm{F}=(\vec{f}_i)_{i=2}^n$ of applied body forces whose tangential component is zero and whose transverse component is directed downwards, as well as a varying number $n$ of vertices.

\begin{figure}
\includegraphics[width = 0.8\textwidth]{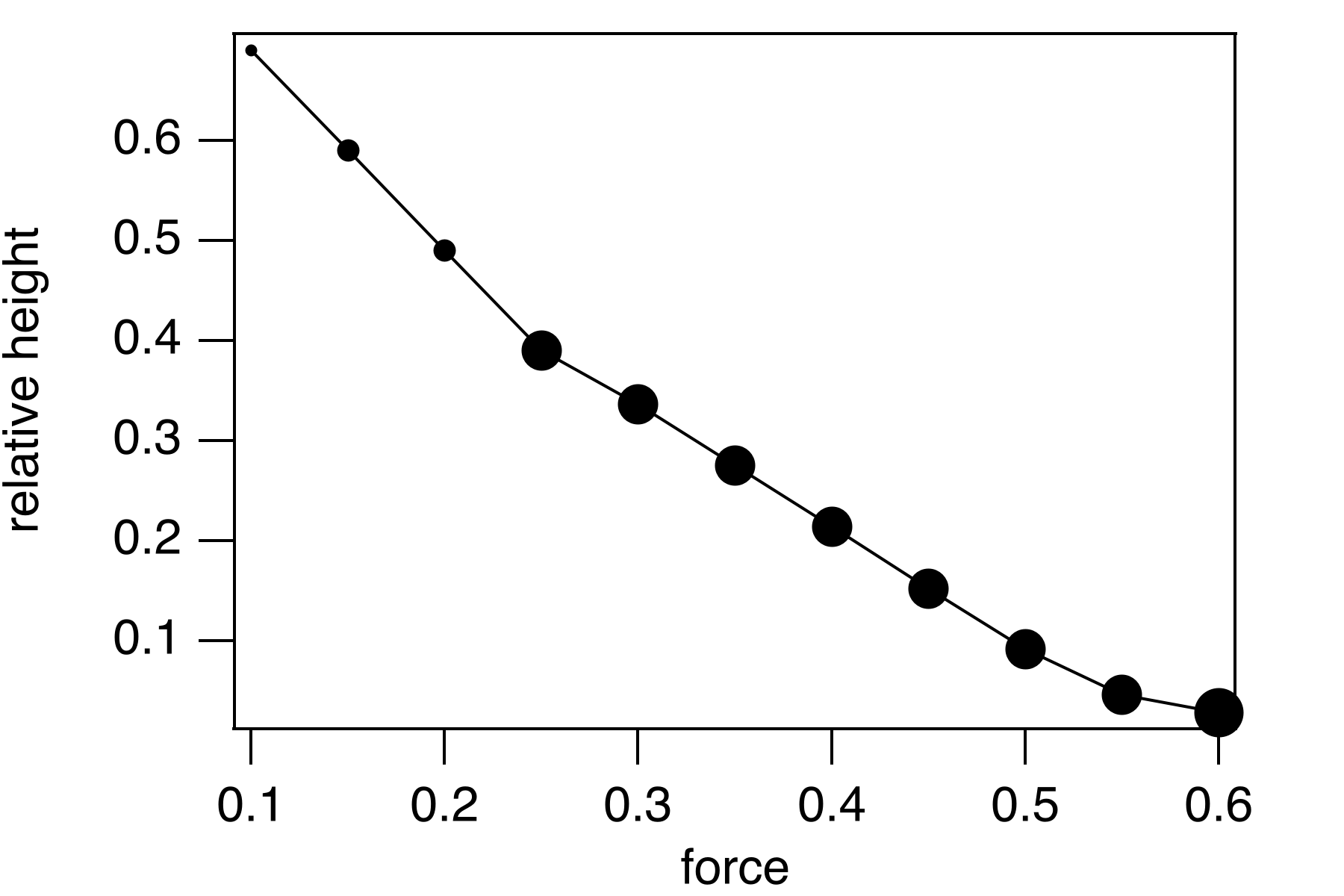}
\caption{Changes in height above the base line, relative to the initial value as the force per vertex is increased. Different marker sizes are associated with different numbers of contact vertices.  }
\label{hvsf}
\end{figure}

\begin{figure}
\includegraphics[width = 0.8\textwidth]{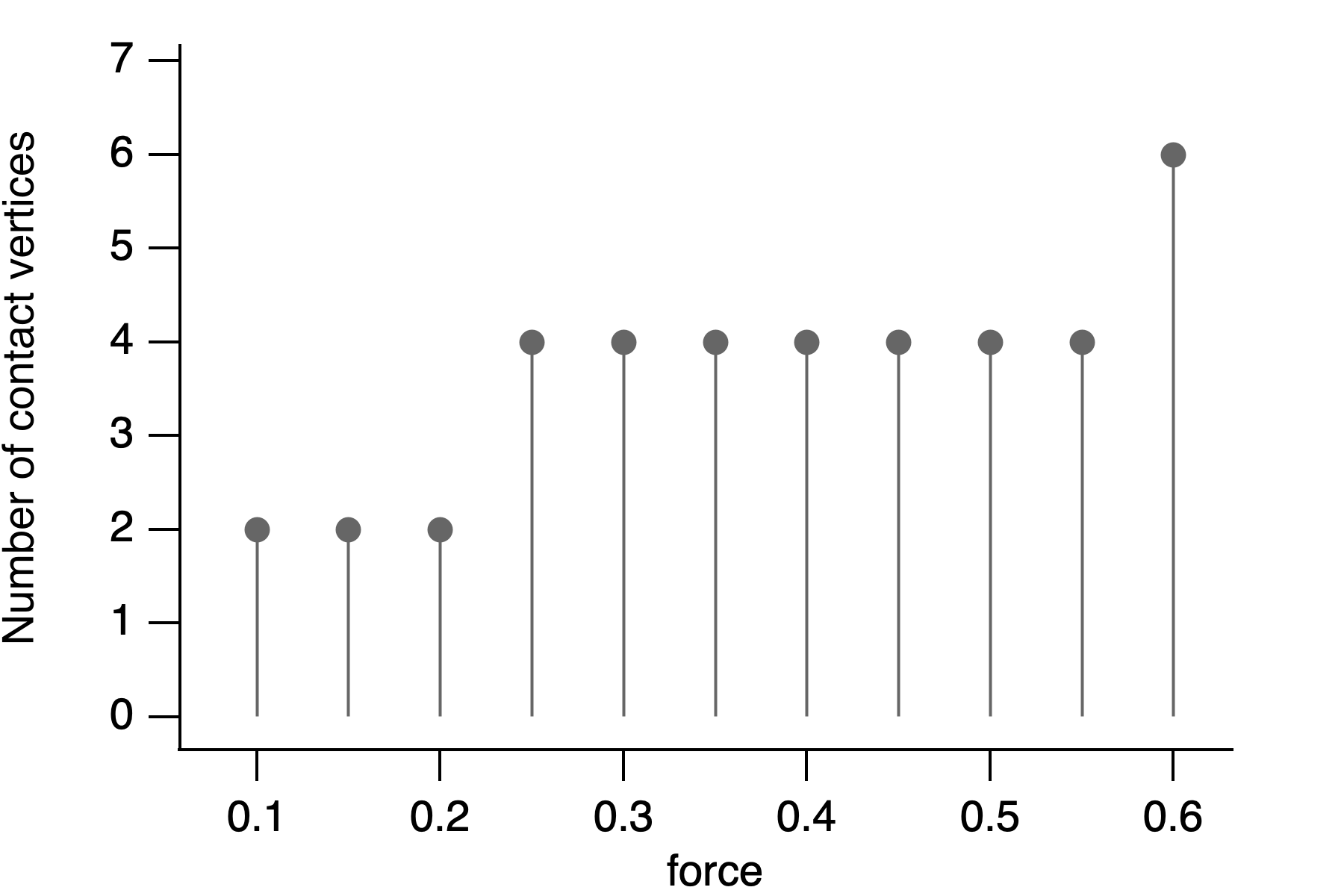}
\caption{Number of contacts vs force for the decagon. Maximum force in this example corresponds to an almost completely flattened structure.}
\label{ncts}
\end{figure}

Figure~\ref{geometry}b shows the shape of a regular decagon under directional pressure for the values $k=\kappa=1.0$, and $f = 0.25$ in reduced units. For this specific load magnitude, half of the vertices are on the surface and half are free. The free vertices align with a good approximation on a circle. The radius of the circle that best approximates the polygonal segment formed of free vertices is 1.9 times the radius of the circle in which the initial regular polygon was inscribed. The top of the deformed polygon is at 0.41 from the initial, unperturbed value. We can now predict how the apparent height above the surface changes as a function of the AFM imaging force. The apparent height is an accurately measurable quantity. Figure~\ref{hvsf} shows the result of such a prediction for the same decagon. The apparent height above the surface decreases smoothly with the applied pressure. This result is consistent with experimental observations obtained from AFM on the brome mosaic virus \cite{Zeng2021}. However, breaks in the slope due to discrete changes in the number of sides laying flat on the base line are noticeable.
While the apparent height varies continuously, the number of surface contacts and hence the surface adhesion does not, Figure~\ref{ncts}. Considering that the imaged structure is the result of the balancing act between adhesive, pressure, and elastic forces, the latter observation predicts that the number of contacts could be discreetly varied by imaging at various probe pressures. Since adhesion is usually strong and irreversible, at least on hydrophobic surfaces, by imaging at low pressures after a higher pressure scan should allow to deduce the contact area from the height and shape of the free surface as the particle relaxes. In other words, for a relaxed polygon (no pressure) with a given partial contact length the shape and the height above the base line should vary discretely, as a function of the number of contacts with the baseline. 

By imaging at high force, and then at low force, and comparing with a model based on the one presented here, one should be able, for the first time, to determine the virus-surface contact area that corresponds to a certain applied force. This is why, we have considered useful to extend our numerical approach and calculate the shape of the polygon, as a function of the number of contact vertices (and thus of the initial imaging force)) once the external force has been turned off, Figure~\ref{geometry}c.

\section{Numerical experiments. Part~{II}: Search for a new equilibrium position for the polygon and correlation with the continuous model}
\label{Sec:4}

In the second experiment, we \emph{assume} that the vertices of the polygon which engage contact with the surface as a result of the application of the applied body force under consideration (cf. section~\ref{Sec:3}) remain in contact with the surface when the applied body force ceases to act. This corresponds to the case where the AFM cantilever pushing the ring downwards is lifted and a certain adhesion force acts on the part of the elastic structure that engaged contact with the flat surface.

As a result, the vertices which do not engage contact with the surface tend to re-organize themselves in a way that minimizes the elastic energy; the points in contact with the surface are only allowed to slide horizontally on the surface. The latter assumption is based on the observation that the activation energy for surface diffusion of adsorbed molecules is usually at least an order of magnitude smaller than the adsorption energy \cite{Somorjai1994}. Thus lateral sliding is much more likely to occur than desorption.

Because of its accurate measurability and practical interest we study the root-mean-square fitted average radius of the free sector (vertices that do not engage with the surface) after the equilibrium restoration of the elastic polygon as a function of the number of contact points with the base line (or, equivalent, as a function of the initial compression described in the previous section). We then verify the genuineness of the proposed discrete model via a set of numerical experiments, where the number of vertices vary. The expectation is that the relaxed shape will tend towards a circular sector that minimizes the line tension.

Table~\ref{rvsn} shows the dependence of the r.m.s. average radius of curvature of the free portion of the decagon on the number of base line contacts. The  changes in the average curvature of the polygon as the number of contact vertices is increased in this example are small, but it is reasonable to assume that they will be measurable if $\sim 10^2$ data points collected during topographic mapping will be used for fitting. As mentioned before, the spatial resolution in the normal direction to the base line is better than 0.5 nm, while the radius of a virus particle is typically 30-100 nm.
\begin{table}
  \centering
  \begin{tabular}{ l | c | c | r }
    \hline
   Number of base line vertices & 3 & 5 & 7 \\ \hline
  R/R$_0$ & 0.99 & 1.05 & 1.04 \\ 
    \hline
  \end{tabular}
   \caption{Ratio of upper polygon average radius (R) and initial, undeformed radius (R$_0$) as a function of the number of baseline vertices.}
  \label{rvsn}
\end{table}

\begin{figure}
\includegraphics[width = \textwidth]{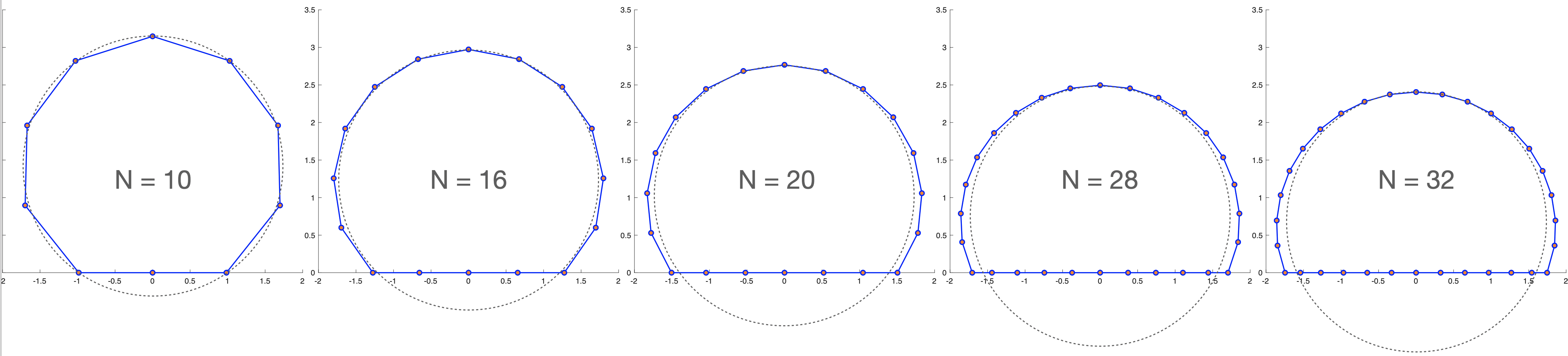}
\caption{Increasing the number of total vertices while keeping all other parameters constant leads to a shape that does not change anymore at higher number of vertices.The dotted line circle corresponds to the circle circumscribing the initial regular polygon and is provided as a reference. }
\label{Nvar}
\end{figure}

Finally, we consider the case where the same total vertical body force is applied, and we let the number of vertices increase. The final polygonal shape should converge to the one expected from a continuum approach, with a circular free perimeter, a fixed contact angle, and height at the top. Figure~\ref{Nvar} shows the result of such a comparison.
Beyond $N \approx 25$ the shape does not change anymore. At this point,  a continuum description will be adequate. The dotted line circle in Figure~\ref{Nvar} corresponds to the circle that circumscribes the initial, regular polyhedron, added to the polygonal shape in a way that superposes the apex point of the circle and the polygon and preserves the left-right symmetry. Below the continuum threshold there are some variations of the polygonal shape with respect to the circumscribing circle, but they are rather minor in the vicinity of the apex point.  However,  the apex point height varies significantly below the continuum threshold. The practical outcome of this observation is that the AFM measurement of the top part of the virus, which is very accurate, could be used in fitting it at every map point with a spherical shell. It is not unreasonable to hypothesize that, from the radius and the center of curvature for that shell one could deduce the elastic constants of the virus particle at an unprecedented level of accuracy. 

\section{Conclusion}
In conclusion, we present a simple model for virus deformation under AFM imaging conditions. The model includes discrete units and deformation modes, adhesion to a substrate, and unidirectional loading forces. We show that numerical solutions based on the model provide a solution uniquely associated to a set of parameters. As a direct application of our analysis, we propose a novel method of using the AFM topographic mapping under different imaging force loads that should lead to determinations of the contact area and elastic constants. The advantages of this approach with respect to state-of-the art elastic deformation studies will be: 1) high throughput, because many viruses can be imaged at the same time, and 2) accuracy, because the resolution of the topographic mapping is very high (sub-nm) in a direction normal to the substrate, even at physiological conditions. 

\section*{Acknowledgments}

We gratefully acknowledge support by the Army Research Office, under award W911NF-17-1-0329 and by the National Science Foundation, under award CBET 1803440 and 1808027 (to BD).  We acknowledge the Center for Bioanalytical Metrology (CBM), an NSF Industry-University Cooperative Research Center, for providing funding under grant NSF IIP 1916645.

This work was partly supported by the Research Fund of Indiana University.

\bibliographystyle{abbrvnat} 
\bibliography{AppAnRef.bib}	

\begin{thebibliography}{46}
\providecommand{\natexlab}[1]{#1}
\providecommand{\url}[1]{\texttt{#1}}
\expandafter\ifx\csname urlstyle\endcsname\relax
  \providecommand{\doi}[1]{doi: #1}\else
  \providecommand{\doi}{doi: \begingroup \urlstyle{rm}\Url}\fi

\bibitem[Aggarwal et~al.(2012)Aggarwal, Rudnick, Bruinsma, and
  Klug]{Aggarwal2012}
A.~Aggarwal, J.~Rudnick, R.~F. Bruinsma, and W.~S. Klug.
\newblock {Elasticity Theory of Macromolecular Aggregates}.
\newblock \emph{Phys. Rev. Lett.}, 109\penalty0 (14):\penalty0 148102, oct
  2012.
\newblock ISSN 0031-9007.
\newblock \doi{10.1103/PhysRevLett.109.148102}.
\newblock URL \url{http://link.aps.org/doi/10.1103/PhysRevLett.109.148102}.

\bibitem[Arkhipov et~al.(2009)Arkhipov, Roos, Wuite, and
  Schulten]{Arkhipov:2009kt}
A.~Arkhipov, W.~H. Roos, G.~J.~L. Wuite, and K.~Schulten.
\newblock {Elucidating the Mechanism behind Irreversible Deformation of Viral
  Capsids}.
\newblock \emph{Biophys. J.}, 97\penalty0 (7):\penalty0 2061--2069, oct 2009.
\newblock ISSN 1542-0086.
\newblock \doi{10.1016/j.bpj.2009.07.039}.
\newblock URL
  \url{http://www.pubmedcentral.nih.gov/articlerender.fcgi?artid=2756377{\&}tool=pmcentrez{\&}rendertype=abstract}.

\bibitem[Baclayon et~al.({2010})Baclayon, Wuite, and Roos]{Baclayon2010}
M.~Baclayon, G.~J.~L. Wuite, and W.~H. Roos.
\newblock {Imaging and manipulation of single viruses by atomic force
  microscopy}.
\newblock \emph{{SOFT MATTER}}, {6}\penalty0 ({21}):\penalty0 {5273--5285},
  {2010}.
\newblock ISSN {1744-683X}.
\newblock \doi{{10.1039/b923992h}}.

\bibitem[Bruinsma and Klug(2015)]{Bruinsma2015}
R.~F. Bruinsma and W.~S. Klug.
\newblock {Physics of Viral Shells}.
\newblock \emph{Annu. Rev. Condens. Matter Phys.}, 6\penalty0 (1):\penalty0
  245--268, mar 2015.
\newblock ISSN 1947-5454.
\newblock \doi{10.1146/annurev-conmatphys-031214-014325}.
\newblock URL
  \url{http://apps.webofknowledge.com/full{\_}record.do?product=UA{\&}search{\_}mode=CitingArticles{\&}qid=2{\&}SID=4BT2Hi2VzgsDBDiKJTy{\&}page=1{\&}doc=7}.

\bibitem[Calo et~al.(2021)Calo, Eleta-Lopez, Ondarcuhu, Verdaguer, and
  Bittner]{Calo2021}
A.~Calo, A.~Eleta-Lopez, T.~Ondarcuhu, A.~Verdaguer, and A.~M. Bittner.
\newblock {Nanoscale wetting of single viruses}.
\newblock \emph{Molecules}, 26\penalty0 (17):\penalty0 5184, aug 2021.
\newblock ISSN 14203049.
\newblock \doi{10.3390/molecules26175184}.
\newblock URL \url{https://www.mdpi.com/1420-3049/26/17/5184}.

\bibitem[Cartagena et~al.(2013)Cartagena, Hernando-P{\'{e}}rez, Carrascosa,
  de~Pablo, and Raman]{Cartagena2013}
A.~Cartagena, M.~Hernando-P{\'{e}}rez, J.~L. Carrascosa, P.~J. de~Pablo, and
  A.~Raman.
\newblock {Mapping in vitro local material properties of intact and disrupted
  virions at high resolution using multi-harmonic atomic force microscopy}.
\newblock \emph{Nanoscale}, 5\penalty0 (11):\penalty0 4729, 2013.
\newblock ISSN 2040-3364.
\newblock \doi{10.1039/c3nr34088k}.
\newblock URL \url{http://xlink.rsc.org/?DOI=c3nr34088k}.

\bibitem[Castellanos et~al.(2012)Castellanos, P{\'{e}}rez, Carrasco,
  Hernando-P{\'{e}}rez, G{\'{o}}mez-Herrero, de~Pablo, and
  Mateu]{Castellanos2012a}
M.~Castellanos, R.~P{\'{e}}rez, C.~Carrasco, M.~Hernando-P{\'{e}}rez,
  J.~G{\'{o}}mez-Herrero, P.~J. de~Pablo, and M.~G. Mateu.
\newblock {Mechanical elasticity as a physical signature of conformational
  dynamics in a virus particle.}
\newblock \emph{Proc. Natl. Acad. Sci. U. S. A.}, 109\penalty0 (30):\penalty0
  12028--33, jul 2012.
\newblock ISSN 1091-6490.
\newblock \doi{10.1073/pnas.1207437109}.
\newblock URL \url{http://www.ncbi.nlm.nih.gov/pubmed/22797893
  http://www.pubmedcentral.nih.gov/articlerender.fcgi?artid=PMC3409779}.

\bibitem[Ciarlet and Piersanti(2019{\natexlab{a}})]{CiaPie2018b}
P.~G. Ciarlet and P.~Piersanti.
\newblock Obstacle problems for {K}oiter's shells.
\newblock \emph{{M}ath. {M}ech. {S}olids}, 24:\penalty0 3061--3079,
  2019{\natexlab{a}}.

\bibitem[Ciarlet and Piersanti(2019{\natexlab{b}})]{CiaPie2018bCR}
P.~G. Ciarlet and P.~Piersanti.
\newblock A confinement problem for a linearly elastic {K}oiter's shell.
\newblock \emph{C.R. Acad. Sci. Paris, S\'{e}r. I}, 357:\penalty0 221--230,
  2019{\natexlab{b}}.

\bibitem[Ciarlet et~al.(2018)Ciarlet, Mardare, and Piersanti]{CiaMarPie2018b}
P.~G. Ciarlet, C.~Mardare, and P.~Piersanti.
\newblock Un probl\`eme de confinement pour une coque membranaire
  lin\'eairement \'elastique de type elliptique.
\newblock \emph{{C}. {R}. {M}ath. {A}cad. {S}ci. {P}aris}, 356\penalty0
  (10):\penalty0 1040--1051, 2018.

\bibitem[Ciarlet et~al.(2019)Ciarlet, Mardare, and Piersanti]{CiaMarPie2018}
P.~G. Ciarlet, C.~Mardare, and P.~Piersanti.
\newblock An obstacle problem for elliptic membrane shells.
\newblock \emph{{M}ath. {M}ech. {S}olids}, 24\penalty0 (5):\penalty0
  1503--1529, 2019.

\bibitem[Dragnea(2019)]{Dragnea2019}
B.~Dragnea.
\newblock {Watching a virus grow}, nov 2019.
\newblock ISSN 10916490.
\newblock URL \url{https://www.pnas.org/content/116/45/22420.short?rss=1}.

\bibitem[Duvaut and Lions(1976)]{DuvLions}
G.~Duvaut and J.-L. Lions.
\newblock \emph{Inequalities in Mechanics and Physics}.
\newblock Springer, Berlin, Heidelberg, 1976.

\bibitem[Ekeland and Temam(1999)]{EkelandTemam1999}
I.~Ekeland and R.~Temam.
\newblock \emph{Convex analysis and variational problems}, volume~28 of
  \emph{Classics in Applied Mathematics}.
\newblock Society for Industrial and Applied Mathematics (SIAM), Philadelphia,
  PA, english edition, 1999.
\newblock Translated from the French.

\bibitem[Gibbons and Klug(2007)]{Gibbons2007}
M.~M. Gibbons and W.~S. Klug.
\newblock {Mechanical modeling of viral capsids}.
\newblock \emph{J. Mater. Sci.}, 42\penalty0 (21):\penalty0 8995--9004, jul
  2007.
\newblock ISSN 0022-2461.
\newblock \doi{10.1007/s10853-007-1741-4}.
\newblock URL
  \url{http://apps.webofknowledge.com/full{\_}record.do?product=UA{\&}search{\_}mode=GeneralSearch{\&}qid=23{\&}SID=4Blr9Lv2MW3Cstzwjc7{\&}page=1{\&}doc=4}.

\bibitem[Greber(2014)]{Greber2014}
U.~F. Greber.
\newblock {How cells tune viral mechanics--insights from biophysical
  measurements of influenza virus.}
\newblock \emph{Biophys. J.}, 106\penalty0 (11):\penalty0 2317--21, jun 2014.
\newblock ISSN 1542-0086.
\newblock \doi{10.1016/j.bpj.2014.04.025}.
\newblock URL
  \url{http://www.pubmedcentral.nih.gov/articlerender.fcgi?artid=4052274{\&}tool=pmcentrez{\&}rendertype=abstract}.

\bibitem[Guerra et~al.(2017)Guerra, Valbuena, Querol-Aud{\'{i}}, Silva,
  Castellanos, Rodr{\'{i}}guez-Huete, Garriga, Mateu, and
  Verdaguer]{Guerra2017}
P.~Guerra, A.~Valbuena, J.~Querol-Aud{\'{i}}, C.~Silva, M.~Castellanos,
  A.~Rodr{\'{i}}guez-Huete, D.~Garriga, M.~G. Mateu, and N.~Verdaguer.
\newblock Structural basis for biologically relevant mechanical stiffening of a
  virus capsid by cavity-creating or spacefilling mutations.
\newblock \emph{Sci. Rep.}, 7\penalty0 (1):\penalty0 4101, dec 2017.
\newblock ISSN 20452322.
\newblock \doi{10.1038/s41598-017-04345-w}.

\bibitem[Ivanovska et~al.(2004)Ivanovska, de~Pablo, Ibarra, Sgalari,
  MacKintosh, Carrascosa, Schmidt, and Wuite]{Ivanovska2004}
I.~L. Ivanovska, P.~J. de~Pablo, B.~Ibarra, G.~Sgalari, F.~C. MacKintosh, J.~L.
  Carrascosa, C.~F. Schmidt, and G.~J.~L. Wuite.
\newblock {Bacteriophage capsids: tough nanoshells with complex elastic
  properties.}
\newblock \emph{Proc. Natl. Acad. Sci. U. S. A.}, 101\penalty0 (20):\penalty0
  7600--5, may 2004.
\newblock ISSN 0027-8424.
\newblock \doi{10.1073/pnas.0308198101}.
\newblock URL \url{http://www.ncbi.nlm.nih.gov/pubmed/15133147
  http://www.pubmedcentral.nih.gov/articlerender.fcgi?artid=PMC419652}.

\bibitem[Kiselev(2006)]{Kiselev2006}
A.~P. Kiselev.
\newblock \emph{Geometry. Book I. Planimetry}.
\newblock Sumizdat, 2006.

\bibitem[Kononova et~al.(2013)Kononova, Snijder, Brasch, Cornelissen, Dima,
  Marx, Wuite, Roos, and Barsegov]{Kononova2013}
O.~Kononova, J.~Snijder, M.~Brasch, J.~Cornelissen, R.~I. Dima, K.~A. Marx,
  G.~J.~L. Wuite, W.~H. Roos, and V.~Barsegov.
\newblock {Structural transitions and energy landscape for Cowpea Chlorotic
  Mottle Virus capsid mechanics from nanomanipulation in vitro and in silico.}
\newblock \emph{Biophys. J.}, 105\penalty0 (8):\penalty0 1893--903, oct 2013.
\newblock ISSN 1542-0086.
\newblock \doi{10.1016/j.bpj.2013.08.032}.

\bibitem[Krishnamani et~al.(2016)Krishnamani, Globisch, Peter, and
  Deserno]{Krishnamani2016}
V.~Krishnamani, C.~Globisch, C.~Peter, and M.~Deserno.
\newblock {Breaking a virus: Identifying molecular level failure modes of a
  viral capsid by multiscale modeling}.
\newblock \emph{Eur. Phys. J. Spec. Top.}, pages 1--18, jul 2016.
\newblock ISSN 1951-6355.
\newblock \doi{10.1140/epjst/e2016-60141-2}.
\newblock URL \url{http://link.springer.com/10.1140/epjst/e2016-60141-2}.

\bibitem[Kuznetsov and McPherson(2011)]{Kuznetsov2011}
Y.~G. Kuznetsov and A.~McPherson.
\newblock {Atomic force microscopy investigation of viruses.}
\newblock \emph{Methods Mol. Biol.}, 736\penalty0 (4):\penalty0 171--195, apr
  2011.

\bibitem[Mannige and Brooks(2009)]{Mannige2009}
R.~V. Mannige and C.~L. Brooks.
\newblock {Geometric considerations in virus capsid size specificity, auxiliary
  requirements, and buckling.}
\newblock \emph{Proc. Natl. Acad. Sci. U. S. A.}, 106\penalty0 (21):\penalty0
  8531--6, may 2009.
\newblock ISSN 1091-6490.
\newblock \doi{10.1073/pnas.0811517106}.
\newblock URL \url{http://www.ncbi.nlm.nih.gov/pubmed/19439655
  http://www.pubmedcentral.nih.gov/articlerender.fcgi?artid=PMC2688982}.

\bibitem[Mateu(2012)]{Mateu2012}
M.~G. Mateu.
\newblock {Mechanical properties of viruses analyzed by atomic force
  microscopy: A virological perspective}.
\newblock \emph{Virus Res.}, 168\penalty0 (1-2):\penalty0 1--22, sep 2012.
\newblock ISSN 01681702.
\newblock \doi{10.1016/j.virusres.2012.06.008}.
\newblock URL
  \url{http://linkinghub.elsevier.com/retrieve/pii/S0168170212002055}.

\bibitem[Michel et~al.(2006)Michel, Ivanovska, Gibbons, Klug, Knobler, Wuite,
  and Schmidt]{Michel:2006jt}
J.~P. Michel, I.~L. Ivanovska, M.~M. Gibbons, W.~S. Klug, C.~M. Knobler,
  G.~J.~L. Wuite, and C.~F. Schmidt.
\newblock {Nanoindentation studies of full and empty viral capsids and the
  effects of capsid protein mutations on elasticity and strength}.
\newblock \emph{Proc. Natl. Acad. Sci. U. S. A.}, 103\penalty0 (16):\penalty0
  6184--6189, apr 2006.
\newblock ISSN 0027-8424.
\newblock \doi{10.1073/pnas.0601744103}.
\newblock URL \url{http://www.pnas.org/content/103/16/6184.full}.

\bibitem[Pang et~al.(2013)Pang, Hevroni, Kol, Eckert, Tsvitov, Kay, and
  Rousso]{Pang2013a}
H.-B. Pang, L.~Hevroni, N.~Kol, D.~M. Eckert, M.~Tsvitov, M.~S. Kay, and
  I.~Rousso.
\newblock {Virion stiffness regulates immature HIV-1 entry.}
\newblock \emph{Retrovirology}, 10:\penalty0 4, jan 2013.
\newblock ISSN 1742-4690.
\newblock \doi{10.1186/1742-4690-10-4}.

\bibitem[Piersanti(2020)]{Pie2020}
P.~Piersanti.
\newblock A time-dependent obstacle problem in linearised elasticity.
\newblock \emph{{N}onlinear {A}nal.}, 192, 2020.

\bibitem[Piersanti(2021)]{Pie2020-1}
P.~Piersanti.
\newblock On the improved interior regularity of the solution of a fourth order
  elliptic problem modelling the displacement of a linearly elastic shallow
  shell subject to an obstacle.
\newblock \emph{{A}symptot. {A}nal.}, 2021.

\bibitem[Piersanti and Shen(2020)]{PS}
P.~Piersanti and X.~Shen.
\newblock Numerical methods for static shallow shells lying over an obstacle.
\newblock \emph{{N}umer. {A}lgorithms}, pages 623--652, 2020.

\bibitem[Piersanti et~al.(In preparation)Piersanti, White, Dragnea, and
  Temam]{PWDT3D}
P.~Piersanti, K.~White, B.~Dragnea, and R.~Temam.
\newblock A three-dimensional discrete model for for approximating the
  deformation of a viral capsid subjected to lying over a flat surface.
\newblock In preparation.

\bibitem[Prevelige et~al.(1993)Prevelige, Thomas, and KING]{Prevelige:1993db}
P.~E. Prevelige, D.~Thomas, and J.~KING.
\newblock {Nucleation and growth phases in the polymerization of coat and
  scaffolding subunits into icosahedral procapsid shells}.
\newblock \emph{Biophys. J.}, 64\penalty0 (3):\penalty0 824--835, mar 1993.

\bibitem[Quarteroni et~al.(2010)Quarteroni, Saleri, and
  Gervasio]{Quarteroni2010}
A.~Quarteroni, F.~Saleri, and P.~Gervasio.
\newblock \emph{Scientific computing with {MATLAB} and {O}ctave}, volume~2 of
  \emph{Texts in Computational Science and Engineering}.
\newblock Springer-Verlag, Berlin, 2010.

\bibitem[Roos(2011)]{Roos2011}
W.~H. Roos.
\newblock {How to perform a nanoindentation experiment on a virus.}
\newblock \emph{Methods Mol. Biol.}, 783:\penalty0 251--64, jan 2011.
\newblock ISSN 1940-6029.
\newblock \doi{10.1007/978-1-61779-282-3_14}.
\newblock URL
  \url{http://apps.webofknowledge.com/full{\_}record.do?product=UA{\&}search{\_}mode=CitingArticles{\&}qid=17{\&}SID=3DdURX6OKNjHSHDYfw8{\&}page=6{\&}doc=59
  http://apps.webofknowledge.com/full{\_}record.do?product=UA{\&}search{\_}mode=GeneralSearch{\&}qid=18{\&}SID=1CKDIiF3nWltHCbh1Rl{\&}page=5{\&}doc=5}.

\bibitem[Sloan(2003)]{Sloan2003}
E.~D. Sloan.
\newblock {Fundamental principles and applications of natural gas hydrates}.
\newblock 426\penalty0 (6964):\penalty0 353--359, nov 2003.
\newblock \doi{10.1038/nature02135}.
\newblock URL \url{http://www.nature.com/articles/nature02135}.

\bibitem[Snijder et~al.(2013)Snijder, Reddy, May, Roos, Nemerow, and
  Wuite]{Snijder2013a}
J.~Snijder, V.~S. Reddy, E.~R. May, W.~H. Roos, G.~R. Nemerow, and G.~J.~L.
  Wuite.
\newblock {Integrin and Defensin Modulate the Mechanical Properties of
  Adenovirus}.
\newblock \emph{J. Virol.}, 87\penalty0 (5):\penalty0 2756--2766, mar 2013.
\newblock ISSN 0022-538X.
\newblock \doi{10.1128/JVI.02516-12}.
\newblock URL \url{http://jvi.asm.org/cgi/doi/10.1128/JVI.02516-12}.

\bibitem[Somorjai(1994)]{Somorjai1994}
G.~A. Somorjai.
\newblock \emph{{Introduction to Surface Chemistry and Catalysis}}.
\newblock John Wiley {\&} Sons, New York, 1994.

\bibitem[Strang(1980)]{Strang1980}
G.~Strang.
\newblock \emph{Linear algebra and its applications}.
\newblock Academic Press [Harcourt Brace Jovanovich, Publishers], New
  York-London, second edition, 1980.

\bibitem[Sun and Yuan(2006)]{SunYuan2006}
W.~Sun and Y.-X. Yuan.
\newblock \emph{Optimization theory and methods}, volume~1 of \emph{Springer
  Optimization and Its Applications}.
\newblock Springer, New York, 2006.
\newblock Nonlinear programming.

\bibitem[Temam and Miranville(2005)]{TemamMiranville2005}
R.~Temam and A.~Miranville.
\newblock \emph{Mathematical modeling in continuum mechanics}.
\newblock Cambridge University Press, Cambridge, second edition, 2005.

\bibitem[Thompson et~al.(2020)Thompson, Cattani, Lee, Ma, Tsvetkova, and
  Dragnea]{Thompson2020}
W.~C. Thompson, A.~J. Cattani, O.~Lee, X.~Ma, I.~B. Tsvetkova, and B.~Dragnea.
\newblock {A Laboratory Model for Virus Particle Nanoindentation}.
\newblock \emph{Biophys.}, 1\penalty0 (2), jan 2020.
\newblock ISSN 2578-6970.
\newblock \doi{10.35459/tbp.2019.000106}.
\newblock URL
  \url{https://meridian.allenpress.com/the-biophysicist/article/doi/10.35459/tbp.2019.000106/436538/A-Laboratory-Model-for-Virus-Particle}.

\bibitem[Vliegenthart and Gompper(2006)]{Vliegenthart2006}
G.~A. Vliegenthart and G.~Gompper.
\newblock {Mechanical deformation of spherical viruses with icosahedral
  symmetry.}
\newblock \emph{Biophys. J.}, 91\penalty0 (3):\penalty0 834--41, aug 2006.
\newblock ISSN 0006-3495.
\newblock \doi{10.1529/biophysj.106.081422}.
\newblock URL \url{http://www.ncbi.nlm.nih.gov/pubmed/16679375
  http://www.pubmedcentral.nih.gov/articlerender.fcgi?artid=PMC1563762}.

\bibitem[Zandi and Reguera(2005)]{Zandi2005}
R.~Zandi and D.~Reguera.
\newblock {Mechanical properties of viral capsids}.
\newblock \emph{Phys. Rev. E}, 72\penalty0 (2):\penalty0 021917, aug 2005.
\newblock ISSN 1539-3755.
\newblock \doi{10.1103/PhysRevE.72.021917}.
\newblock URL \url{http://link.aps.org/doi/10.1103/PhysRevE.72.021917}.

\bibitem[Zeng et~al.(2017)Zeng, Hernando-P{\'{e}}rez, Dragnea, Ma, van~der
  Schoot, and Zandi]{Zeng2017a}
C.~Zeng, M.~Hernando-P{\'{e}}rez, B.~Dragnea, X.~Ma, P.~van~der Schoot, and
  R.~Zandi.
\newblock {Contact Mechanics of a Small Icosahedral Virus}.
\newblock \emph{Phys. Rev. Lett.}, 119\penalty0 (3):\penalty0 038102, jul 2017.
\newblock ISSN 0031-9007.
\newblock \doi{10.1103/PhysRevLett.119.038102}.
\newblock URL \url{http://link.aps.org/doi/10.1103/PhysRevLett.119.038102}.

\bibitem[Zeng et~al.(2021)Zeng, Scott, Malyutin, Zandi, {Der Schoot}, and
  Dragnea]{Zeng2021}
C.~Zeng, L.~Scott, A.~Malyutin, R.~Zandi, P.~V. {Der Schoot}, and B.~Dragnea.
\newblock {Virus mechanics under molecular crowding}.
\newblock \emph{J. Phys. Chem. B}, 125\penalty0 (7):\penalty0 1790--1798, feb
  2021.
\newblock ISSN 15205207.
\newblock \doi{10.1021/acs.jpcb.0c10947}.
\newblock URL \url{https://pubs.acs.org/doi/10.1021/acs.jpcb.0c10947}.

\bibitem[Zlotnick(2007)]{Zlotnick:2007ji}
A.~Zlotnick.
\newblock {Distinguishing reversible from irreversible virus capsid assembly.}
\newblock \emph{J. Mol. Biol.}, 366\penalty0 (1):\penalty0 14--18, feb 2007.

\bibitem[Zlotnick et~al.(1999)Zlotnick, Johnson, Wingfield, Stahl, and
  Endres]{Zlotnick:1999ww}
A.~Zlotnick, J.~M. Johnson, P.~W. Wingfield, S.~J. Stahl, and D.~Endres.
\newblock {A theoretical model successfully identifies features of hepatitis B
  virus capsid assembly.}
\newblock \emph{Biochemistry}, 38\penalty0 (44):\penalty0 14644--14652, nov
  1999.

\end{thebibliography}

\end{document}